\documentclass{amsart}

\newcommand\inprod[2]{\left<#1,#2\right>}

\usepackage{amssymb}
\usepackage{amsmath}
\usepackage{amsthm}
\newtheorem{Theorem}{Theorem}
\usepackage{mathenv}
\usepackage{amsfonts}
\usepackage{verbatim}

\title{Uniqueness of the 2-universality Criterion}
\author{Scott D. Kominers}\thanks{This research was supported by a fellowship from the Center for Excellence in Education and was conducted at the 2004 Research Science Institute.}
\address{Department of Mathematics, Harvard University\\8520 Burning Tree Road, Bethesda, MD 20817}\email{kominers@fas.harvard.edu}
\date{\today}

\begin{document}
\begin{abstract}Kim, Kim, and Oh gave a minimal criterion for the $2$-universality of positive-definite integer-matrix quadratic forms.  We show that this $2$-universality criterion is unique in the sense of the uniqueness of the Conway-Schneeberger Fifteen Theorem.\end{abstract}
\maketitle
\section{Introduction}
By a \textit{quadratic form} (or just \emph{form}) \emph{of rank $n$} we mean a degree-two homogenous polynomial in $n$ independent variables.   If the  quadratic form $Q$ is given by $Q(x_1,...,x_n)= \sum_{i,j}a_{ij}x_ix_j$ with $a_{ij}=a_{ji}$, then the matrix given by $L=(a_{ij})$ is the \emph{Gram Matrix} of a $\mathbb{Z}$-lattice $L$ equipped with a  symmetric bilinear form $\inprod{\cdot}{\cdot}$ such that $\inprod{L}{L}\subseteq \mathbb{Z}$.  We have immediately from these structures that $Q(\mathbf{x})=\mathbf{x}^TL\mathbf{x}=\inprod{L\mathbf{x}}{\mathbf{x}}$ for $\mathbf{x}\in\mathbb{R}^n$.

For convenience, we use form-theoretic and lattice-theoretic language interchangeably throughout.  A complete introduction to both approaches to quadratic form theory can be found in \cite{O'meara:Lattice}.

We say that a rank-$n$ form $Q$ \textit{represents} an integer $k$ if there is an $\mathbf{x}\in\mathbb{Z}^n$ such that $Q(\mathbf{x})=k$.  More generally, we say that a lattice $L$ represents another lattice $\ell$ if there is a $\mathbb{Z}$-linear, bilinear form-preserving injection $\sigma:\ell \to L$.  A form is called \textit{universal} if it represents all positive integers and is similarly called \textit{$n$-universal} if it represents all positive-definite integer-matrix rank-$n$ quadratic forms. It is clear that a rank-$n$ form $Q$ is universal if and only if it is $1$-universal, as for an integer $k$ $$k = Q(x_1,\ldots, x_n)\iff Q(x_1x,\ldots,x_nx)=kx^2.$$

In 1993, Conway and Schneeberger announced the \emph{Fifteen Theorem}, giving a criterion characterizing the positive-definite integer-matrix quadratic forms which represent all positive integers.   Specifically, they showed that any positive-definite integer-matrix form which represents the set of nine critical numbers $$\mathcal{S}_1=\{1,2,3,5,6,7,10,14,15\}$$ is universal \cite{Conway:universality,Bhargava:Fif}.  Kim, Kim, and Oh \cite{Kim:universal} presented an analogous criterion for $2$-universality which we state in Theorem~\ref{2-crit} of Section~\ref{2-critSec}.

The set $\mathcal{S}_1$ of the Fifteen Theorem is known to be unique.  Indeed, if $\mathcal{S}_1'$ is a set of integers such that a quadratic form is universal if and only if it represents the full set $\mathcal{S}_1'$, then $\mathcal{S}_1\subseteq \mathcal{S}_1'$.  We show an analogous uniqueness result for the 2-universality criterion found by Kim, Kim, and Oh~\cite{Kim:universal}.

\section{Notations and Terminology}
If a $\mathbb{Z}$-lattice $L$ is of the form $L= L_1\oplus L_2$ for sublattices $L_1$, $L_2$ of $L$ and $\inprod{L_1}{L_2}=0$ then we write $L\cong L_1\bot L_2$ and say that $L_1$ and $L_2$ are \textit{orthogonal}.    

We write $\left<a_1,\ldots,a_n\right>$ for the rank-$n$ diagonal form  $$a_1x_1^2+\cdots+a_nx_n^2\cong \left(\begin{array}{ccc}a_1&&\\&\ddots&\\&&a_n\end{array}\right)$$   and denote  by $[a,b,c]$ the rank-$2$ form $$ax^2+2bxy+cy^2\cong \left(\begin{array}{cc}a&b\\b&c\end{array}\right).$$  From the classical reduction theory of quadratic forms, we may assume that the form $[a,b,c]$ is always \emph{Minkowski-reduced} so that $0\leq 2b\leq a\leq c$.

We work with a generalization of the escalation method used by Conway~\cite{Conway:universality} and Bhargava~\cite{Bhargava:Fif}.  Extending the definitions of Bhargava \cite{Bhargava:Fif}, we define a \textit{truant}  of a lattice $L$ to be a lattice not represented by $L$.  An \textit{escalation} of $L$ by a rank-$n$ truant $\ell$ is a lattice $L'$ representing $\ell$ which contains $L$ as a sublattice with codimension at most $n$.  

If $\mathcal{S}$ is a set of rank-$n$ forms such that all escalations by elements in $\mathcal{S}$ eventually produce lattices which are $n$-universal, then every lattice which represents all of $\mathcal{S}$ must contain an $n$-universal sublattice and thus is itself $n$-universal (see \cite{Conway:universality, Bhargava:Fif, Kim:finite}).  We call any such $\mathcal{S}$ an \emph{$n$-criterion set}.  Thus, for example, the set $\mathcal{S}_1$ found by Conway \cite{Conway:universality} naturally gives the $1$-criterion set $$\{x^2,2x^2,3x^2,5x^2,6x^2,7x^2,10x^2,14x^2,15x^2\}.$$

\section{Uniqueness of the 2-criterion Set}
\label{2-critSec}
Kim, Kim, and Oh found the following $2$-criterion set in \cite{Kim:universal}:
\begin{Theorem}[Kim, Kim, and Oh] \label{2-crit}A $2$-criterion set is given by $$\mathcal{S}_2:=\{\left<1,1\right>,\ \left<2,3\right>,\ \left<3,3\right>,\ [2,1,2],\ [2,1,3],\ [2,1,4]\}.$$\end{Theorem}

More can be said about this criterion: the set $\mathcal{S}_2$ is a \textit{minimal} 2-criterion set, in the sense that for every form $\ell\in \mathcal{S}_2$ there is some rank-4 form which represents all of $\mathcal{S}_2$ but $\ell$ (see \cite{Kim:universal}).  We now strengthen this result, showing that $\mathcal{S}_2$ is the \textit{unique} minimal 2-criterion set.
\begin{Theorem}\label{unique}The set of forms $\mathcal{S}_2$ given in Theorem~\ref{2-crit} is the unique minimal 2-criterion set---that is, every $2$-criterion set must contain $\mathcal{S}_2$ as a subset.
\end{Theorem}

\begin{proof}
Throughout, $\mathcal{T}$ denotes a finite set of rank-2 forms not containing some form $\ell\in \mathcal{S}_2$.  It suffices to show that for any such $\mathcal{T}$ there is some lattice with truant $\ell$ which represents all of $\mathcal{T}$, since we know from Theorem \ref{2-crit} that $\mathcal{S}_2$ is a $2$-criterion set.

If $\left<1,1\right>\not\in \mathcal{T}$ then we may write (by Minkowski reduction) $$\mathcal{T}=\{\left<1,c_1\right>,\ldots,\left<1,c_k\right>, L_1,\ldots, L_{k'}\},$$ where $c_i>1$ for all $1\leq i\leq k$ and the first minimum of $L_i$ is also larger than $1$ for each $1\leq i\leq k'$.  Then, the lattice $$\left<1,c_1,\ldots,c_k\right>\bot L_1\bot \ldots \bot L_{k'}$$ represents all of $\mathcal{T}$ but has truant $\left<1,1\right>$.  We have therefore shown that any $2$-criterion set must contain $\left<1,1\right>$.

Now, if $\left<2,3\right>\not\in \mathcal{T}$ then we may express \begin{align*}\{\left<a_1,c_1\right>,\ldots,\left<a_k,c_k\right> \}&:=\left\{\left<a,c\right>\in\mathcal{T}\mid a\in \{1,2,3\}, c>3\right\},\\
\{[d_1,1,e_1],\ldots,[d_{k'},1,e_{k'}]\}&:=\left\{[d,1,e]\in\mathcal{T}\mid d\in \{2,3\}, e>4\right\},\\
\{L_1,\ldots, L_{k''}\}&:=\left\{[p,q,r]\in\mathcal{T}\mid 3<p\leq r\right\}.
\end{align*}Then, the lattice $$\left<1,1,c_1,\ldots,c_k\right>\bot [2,1,e_1]\bot\ldots\bot[2,1,e_{k'}]\bot L_1\bot \ldots \bot L_{k''}$$ represents all of $\mathcal{T}$ but has truant $\left<2,3\right>$, whence every $2$-criterion set must contain $\left<2,3\right>$.  An analogous argument shows that every $2$-criterion set must also contain  $\left<3,3\right>$.

Likewise, if $[2,1,e_*]\not\in\mathcal{T}$ for some $e_*\in\{2,3,4\}$ then we consider the sets \begin{align*}\{\left<a_1,c_1\right>,\ldots,\left<a_k,c_k\right>\}&:=\left\{\left<a,c\right>\in\mathcal{T}\mid a\geq 2, c>e_*\right\},\\
\{[d_1,1,e_1],\ldots,[d_{k'},1,e_{k'}]\}&:=\left\{[d,1,e]\in\mathcal{T}\mid d\geq 2, e>e_*\right\},\\
\{L_1,\ldots, L_{k''}\}&:=\left\{[p,q,r]\in\mathcal{T}\mid 3<p\leq r\right\}.
\end{align*}As the rank-$e_*$ form $\left<1,\ldots,1\right>$ represents $[2,1,e]$ for all $1<e<e_*$, we observe that the lattice $$\underbrace{\left<1,\ldots,1\right>}_{e_* \text{ times}}\bot \left<a_1,c_1,\ldots, a_k,c_k\right>\bot[d_1,1,e_1]\bot \ldots \bot [d_{k'},1,e_{k'}]\bot L_1\bot \ldots \bot L_{k''} $$
represents all of $\mathcal{T}$ but does not represent $[2,1,e_*]$.  We therefore see that every $2$-criterion set must contain  $[2,1,e_*]$ for each $e_*\in\{2,3,4\}$.

Since we shown that every $2$-criterion set must include each $\ell\in\mathcal{S}_2$, we have proven the theorem.
\end{proof}

\section*{Acknowledgements}
The author is grateful to Aaron M. Tievsky, who supervised the research and gave comments on several earlier versions of this article.  He also thanks Zachary Abel, Noam D. Elkies, Ellen Dickstein Kominers, Paul M. Kominers, William Kominers, Kai-Wen Lan, Christopher C. Mihelich, Chung Pang Mok, Anatoly Preygel, Antoni Rangachev, John H. Rickert,  Evgenia Sendova, Linda B. Westrick, Susan Schwartz Wildstrom, and an anonymous referee for their helpful comments and suggestions on the work and on earlier drafts of this article.

\end{document}